\documentclass{amsart}
\usepackage{amssymb}

\usepackage{graphicx}
\usepackage{amscd}

\newtheorem{theorem}{Theorem}
\theoremstyle{plain}

\newtheorem{corollary}{Corollary}

\newtheorem{definition}{Definition}

\newtheorem{proposition}{Proposition}

\numberwithin{equation}{section}

\input{tcilatex}
\input{diagrams.tex}
\newarrow{Onto}{-}{-}{-}{-}{->>}
\newarrow{Into}{C}{-}{-}{-}{>}
\newarrow{Corr}{<}{-}{-}{-}{>}
\newarrow{Ev}{|}{-}{-}{-}{>}
\begin{document}
\title{Generalized Twisted coHom Objects}
\author{Sergio D. Grillo}
\address{\textit{Centro At\'{o}mico Bariloche and Instituto Balseiro}\\
\textit{\ 8400-S. C. de Bariloche}\\
\textit{\ Argentina}}
\email{sergio@cabtep2.cnea.gov.ar}

\begin{abstract}
A generalization of the concept of twisted internal coHom object in the
category of conic quantum spaces (c.f. \cite{gm}) was outlined in \cite{g}.
The aim of this article is to develope in more detail this generalization.
\end{abstract}

\maketitle

\section{Introduction}

Given a couple $\mathcal{A},\mathcal{B}$ of conic quantum spaces, i.e. $%
\mathcal{A},\mathcal{B}\in \mathrm{CA}$ (the monoidal category of finitely
generated graded algebras \cite{man0}\cite{man1}\cite{gm}), their symmetric
twisted tensor products $\mathcal{A}\circ _{\tau }\mathcal{B}$ \cite{cap}
\cite{gm} can also be seen as (2nd admissible) counital 2-cocycle twisting
of the quantum space $\mathcal{A}\circ \mathcal{B}$ \cite{g}. Then, in the
same sense as in \cite{gm}, we can study, instead of the maps $\mathcal{A}%
\rightarrow \mathcal{H}\circ \mathcal{B}$ (which define a comma category
whose initial objects are the proper internal coHom objects of $\mathrm{CA}$%
), certain subclasses of arrows $\mathcal{A}\rightarrow \left( \mathcal{H}%
\circ \mathcal{B}\right) _{\omega }$, where $\omega $ is a counital
2-cocycle defining a twist transformation on $\mathcal{H}\circ \mathcal{B}$.
The aim of this paper is to show that, in certain circumstances, these
classes give rise for each pair $\mathcal{A},\mathcal{B}\in \mathrm{CA}$ to
a category $\Omega ^{\mathcal{A},\mathcal{B}}$ with initial object, namely $%
\underline{hom}^{\Omega }\left[ \mathcal{B},\mathcal{A}\right] $, such that
the disjoint union $\Omega ^{\cdot }=\bigvee_{\mathcal{A},\mathcal{B}}\Omega
^{\mathcal{A},\mathcal{B}}$ has a semigroupoid structure together with a
related embedding $\Omega ^{\cdot }\hookrightarrow \mathrm{CA}$ that
preserves the involved (partial) products. Consequently, $\left( \mathcal{B},%
\mathcal{A}\right) \mapsto \underline{hom}^{\Omega }\left[ \mathcal{B},%
\mathcal{A}\right] $ defines an $\mathrm{CA}$-cobased category with an
additional notion of evaluation given by arrows
\begin{equation*}
\mathcal{A}\rightarrow \left( \underline{hom}^{\Omega }\left[ \mathcal{B},%
\mathcal{A}\right] \circ \mathcal{B}\right) _{\omega }.
\end{equation*}
The categories $\Upsilon ^{\cdot }$, obtained in \cite{gm}, are particular
cases of the categories $\Omega ^{\cdot }$. In this way we generalize the
idea of twisted coHom objects in the more general framework of twisting of
quantum spaces. This setting enable us, in turn, to a better understanding
of the results obtained in the mentioned paper.

This article is based on the contents of \cite{gm} and \cite{g} and
references therein, thus we shall frequently refer the reader to them.
Notation and terminology also follow those papers.

\section{The categories $\Omega ^{\cdot }$}

In order to built up the categories $\Omega ^{\mathcal{A},\mathcal{B}}$, let
us first make a couple of observations.

1. Consider on the category $\mathrm{GrVct}$ (of graded vector spaces) the
monoid
\begin{equation*}
\mathbf{V}\circ \mathbf{W}=\bigoplus_{n\in \mathbb{N}_{0}}\left( \mathbf{V}%
_{n}\otimes \mathbf{W}_{n}\right) ,
\end{equation*}
for $\mathbf{V}=\bigoplus_{n\in \mathbb{N}_{0}}\mathbf{V}_{n}$ and $\mathbf{W%
}=\bigoplus_{n\in \mathbb{N}_{0}}\mathbf{W}_{n}$. The arrows $\mathbf{V}%
\rightarrow \mathbf{W}$ in $\mathrm{GrVct}$ are homogeneous linear maps;
i.e. its restrictions to each $\mathbf{V}_{n}$ define maps $\mathbf{V}%
_{n}\rightarrow \mathbf{W}_{n}$. It is clear that the forgetful functor $%
\frak{H}:\mathrm{CA}\hookrightarrow \mathrm{GrVct}$ turns into a monoidal
one, and $\frak{H}\left( \mathcal{A}_{\psi }\right) =\frak{H}\left( \mathcal{%
A}\right) =\mathbf{A}=\bigoplus_{n\in \mathbb{N}_{0}}\mathbf{A}_{n}$ holds
for every twist transformation $\psi \in \frak{Z}^{2}\left[ \mathbf{A}_{1}%
\right] $.

2. Let us construct the comma categories $\left( \frak{H}\left( \mathcal{A}%
\right) \downarrow \frak{H\,}\left( \mathrm{CA}\circ \mathcal{B}\right)
\right) $, where the functor $\frak{H\,}\left( \mathrm{CA}\circ \mathcal{B}%
\right) $ is the composition of $\mathrm{CA}\circ \mathcal{B}$ and $\frak{H}$%
. Its objects are pairs $\left\langle \varphi ,\mathcal{H}\right\rangle $
where $\mathcal{H}\in \mathrm{CA}$ and $\varphi $ is an arrow in $\mathrm{%
GrVct}$,
\begin{equation*}
\varphi :\frak{H}\left( \mathcal{A}\right) \rightarrow \frak{H}\left(
\mathcal{H}\circ \mathcal{B}\right) =\frak{H}\left( \mathcal{H}\right) \circ
\frak{H}\left( \mathcal{B}\right) .
\end{equation*}
To every $\left\langle \varphi ,\mathcal{H}\right\rangle \in \left( \frak{H}%
\left( \mathcal{A}\right) \downarrow \frak{H\,}\left( \mathrm{CA}\circ
\mathcal{B}\right) \right) $, the surjection $\pi ^{\varphi }:\mathbf{B}%
_{1}^{\ast }\otimes \mathbf{A}_{1}\twoheadrightarrow \mathbf{H}_{1}^{\varphi
}:b^{j}\otimes a_{i}\mapsto h_{i}^{j}$ can be related, being $h_{i}^{j}$ the
elements of $\mathbf{H}_{1}$ defining the restriction of $\varphi $ to $%
\mathbf{A}_{1}$, i.e. $\varphi \left( a_{i}\right) =h_{i}^{j}\otimes b_{j}$.
This linear surjection gives rise to a functor
\begin{equation}
\begin{array}{c}
\frak{F}:\left( \frak{H}\left( \mathcal{A}\right) \downarrow \frak{H\,}%
\left( \mathrm{CA}\circ \mathcal{B}\right) \right) \rightarrow \mathrm{CA},
\\
\\
\left\langle \varphi ,\mathcal{H}\right\rangle \mapsto \mathcal{H}^{\varphi
}=\left( \mathbf{H}_{1}^{\varphi },\mathbf{H}^{\varphi }\right)
;\;\;\;\alpha \mapsto \left. \alpha \right| _{\mathbf{H}^{\varphi }},
\end{array}
\label{f}
\end{equation}
where $\mathbf{H}^{\varphi }$ is the subalgebra of $\mathbf{H}$ generated by
$\mathbf{H}_{1}^{\varphi }$ (an analogous functor is used in \cite{gm} to
built up the categories $\Upsilon ^{\mathcal{A},\mathcal{B}}$).

\bigskip

Using last functor we shall construct each $\Omega ^{\mathcal{A},\mathcal{B}}
$ as a full subcategory of the corresponding comma category $\left( \frak{H}%
\left( \mathcal{A}\right) \downarrow \frak{H\,}\left( \mathrm{CA}\circ
\mathcal{B}\right) \right) $. Given a couple of conic quantum spaces $%
\mathcal{A}$ and $\mathcal{B}$, consider a counital element
\begin{equation}
\omega :\left( \left[ \mathbf{B}_{1}^{\ast }\otimes \mathbf{A}_{1}\otimes
\mathbf{B}_{1}\right] ^{\otimes }\right) ^{\otimes 2}\backsimeq \left( \left[
\mathbf{B}_{1}^{\ast }\otimes \mathbf{A}_{1}\otimes \mathbf{B}_{1}\right]
^{\otimes }\right) ^{\otimes 2}  \label{lb}
\end{equation}
of $\frak{Z}^{2}\left[ \mathbf{B}_{1}^{\ast }\otimes \mathbf{A}_{1}\otimes
\mathbf{B}_{1}\right] $. Eventually, for $\left\langle \varphi ,\mathcal{H}%
\right\rangle \in \left( \frak{H}\left( \mathcal{A}\right) \downarrow \frak{%
H\,}\left( \mathrm{CA}\circ \mathcal{B}\right) \right) $, we can translate $%
\omega $ to $\left[ \mathbf{H}_{1}^{\varphi }\otimes \mathbf{B}_{1}\right]
^{\otimes }$ through $\pi ^{\varphi }:\mathbf{B}_{1}^{\ast }\otimes \mathbf{A%
}_{1}\twoheadrightarrow \mathbf{H}_{1}^{\varphi }$ in such a way that the
diagram
\begin{equation*}
\begin{diagram} \QTR{bf}{B}_{1}\otimes\QTR{bf}{B}_{1}^{\ast }\otimes
\QTR{bf}{A}_{1} & \rOnto &\QTR{bf}{B}_{1}\otimes \QTR{bf}{H}_{1}^{\varphi
}\\ \dCorr~{\backsimeq } & & \dTo \\
\QTR{bf}{B}_{1}\otimes\QTR{bf}{B}_{1}^{\ast }\otimes \QTR{bf}{A}_{1} &
\rOnto &\QTR{bf}{B}_{1}\otimes \QTR{bf}{H}_{1}^{\varphi }\\ \end{diagram}
\end{equation*}
be commutative, defining in this way a counital 2-cochain in $\frak{Z}^{2}%
\left[ \mathbf{H}_{1}^{\varphi }\otimes \mathbf{B}_{1}\right] $. Last
affirmation lies on the results given in \textbf{Prop. 4 }of \cite{g},
applied to the injection $\mathbf{H}_{1}^{\varphi }\hookrightarrow \mathbf{B}%
_{1}^{\ast }\otimes \mathbf{A}_{1}$. Then, if the resulting automorphism is
admissible, we can define with it a twist transformation on $\mathcal{H}%
^{\varphi }\circ \mathcal{B}=\frak{F}\left\langle \varphi ,\mathcal{H}%
\right\rangle \circ \mathcal{B}$.

\begin{definition}
For every pair $\mathcal{A},\mathcal{B}\in \mathrm{CA}$ and $\omega \in
\frak{Z}^{2}\left[ \mathbf{B}_{1}^{\ast }\otimes \mathbf{A}_{1}\otimes
\mathbf{B}_{1}\right] $, we define $\Omega ^{\mathcal{A},\mathcal{B}}$ as
the full subcategory of $\left( \frak{H}\left( \mathcal{A}\right) \downarrow
\frak{H\,}\left( \mathrm{CA}\circ \mathcal{B}\right) \right) $ formed out by
diagrams $\left\langle \varphi ,\mathcal{H}\right\rangle $ such that $\omega
$ defines a 2-cocycle twisting on $\frak{F}\left\langle \varphi ,\mathcal{H}%
\right\rangle \circ \mathcal{B}$, and the homogeneous linear map $\varphi $
is a morphism of quantum spaces $\mathcal{A}\rightarrow \left( \frak{F}%
\left\langle \varphi ,\mathcal{H}\right\rangle \circ \mathcal{B}\right)
_{\omega }.\;\;\;\;\blacksquare $
\end{definition}

Given now a collection of cochains $\left\{ \omega _{\mathcal{A},\mathcal{B}%
}\right\} _{\mathcal{A},\mathcal{B}\in \mathrm{CA}}\subset \frak{Z}^{2}\left[
\mathbf{B}_{1}^{\ast }\otimes \mathbf{A}_{1}\otimes \mathbf{B}_{1}\right] $,
we name $\Omega ^{\cdot }$ the disjoint union of the categories $\Omega ^{%
\mathcal{A},\mathcal{B}}$ just defined. Clearly, $\mathrm{CA}^{\circ }$ (see
\cite{gm}, or \S \textbf{1.2} of \cite{g} for a brief review) is a category $%
\Omega ^{\cdot }$ with an associated collection given by identity maps.

Calling $\frak{H}\mathrm{CA}^{\circ }$ the disjoint union of $\left( \frak{H}%
\left( \mathcal{A}\right) \downarrow \frak{H\,}\left( \mathrm{CA}\circ
\mathcal{B}\right) \right) $, it follows that every $\Omega ^{\cdot }$ is a
full subcategory of $\frak{H}\mathrm{CA}^{\circ }$. On the other hand, let
us observe that $\frak{H}\mathrm{CA}^{\circ }$ has a semigroupoid structure
given by the functor
\begin{equation}
\left\langle \varphi ,\mathcal{H}\right\rangle \times \left\langle \chi ,%
\mathcal{G}\right\rangle \mapsto \left\langle \left( I_{H}\circ \chi \right)
\,\varphi ,\mathcal{H}\circ \mathcal{G}\right\rangle ;\;\;\alpha \times
\beta \mapsto \alpha \circ \beta ,  \label{ma}
\end{equation}
and $\mathrm{CA}^{\circ }\subset \frak{H}\mathrm{CA}^{\circ }$ is a\textbf{\
}sub-semigroupoid. In fact, this map is a partial product functor with
domain
\begin{equation*}
\bigvee_{\mathcal{A},\mathcal{B},\mathcal{C}\in \mathrm{CA}}\left( \frak{H}%
\left( \mathcal{A}\right) \downarrow \frak{H}\left( \mathrm{CA}\circ
\mathcal{C}\right) \right) \times \left( \frak{H}\left( \mathcal{C}\right)
\downarrow \frak{H}\left( \mathrm{CA}\circ \mathcal{B}\right) \right)
\end{equation*}
and codomain $\frak{H}\mathrm{CA}^{\circ }$, such that
\begin{equation*}
\left( \frak{H}\left( \mathcal{A}\right) \downarrow \frak{H}\left( \mathrm{CA%
}\circ \mathcal{C}\right) \right) \times \left( \frak{H}\left( \mathcal{C}%
\right) \downarrow \frak{H}\left( \mathrm{CA}\circ \mathcal{B}\right)
\right) \rightarrow \left( \frak{H}\left( \mathcal{A}\right) \downarrow
\frak{H}\left( \mathrm{CA}\circ \mathcal{B}\right) \right)
\end{equation*}
Its associativity comes from that of $\circ $, and the unit elements are
given by the diagrams $\left\langle \ell _{\mathcal{A}},\mathcal{K}%
\right\rangle $, where $\ell _{\mathcal{A}}$ is the homogeneous isomorphism $%
\mathbf{A}\backsimeq \Bbbk \left[ e\right] \otimes \mathbf{A}$, such that $%
a\mapsto e^{n}\otimes a$ if $a\in \mathbf{A}_{n}$. Nevertheless, for a
generic collection $\left\{ \omega _{\mathcal{A},\mathcal{B}}\right\} _{%
\mathcal{A},\mathcal{B}\in \mathrm{CA}}$ of cochains, $\Omega ^{\cdot }$
fails to be a semigroupoid. Furthermore, in the generic case, each $\Omega ^{%
\mathcal{A},\mathcal{B}}$ fails to have initial objects. To address this
problem, we shall consider particular cases.

\section{The semigroupoid structure of $\Omega ^{\cdot }$}

In what follows, all references to sections and theorems correspond to \cite
{g}. Recall the monics (c.f. \S \textbf{2.3.2})
\begin{equation}
\frak{j}:\frak{C}^{\bullet }\left[ \mathbf{B}_{1}\right] ^{!}\times \frak{C}%
^{\bullet }\left[ \mathbf{A}_{1}\right] \times \frak{C}^{\bullet }\left[
\mathbf{B}_{1}\right] \hookrightarrow \frak{C}^{\bullet }\left[ \mathbf{B}%
_{1}^{\ast }\otimes \mathbf{A}_{1}\otimes \mathbf{B}_{1}\right] .  \label{m3}
\end{equation}

\begin{definition}
A collection $\left\{ \omega _{\mathcal{A},\mathcal{B}}\right\} _{\mathcal{A}%
,\mathcal{B}\in \mathrm{CA}}$ is \textbf{factorizable} if there exists
another collection
\begin{equation*}
\left\{ \psi _{\mathcal{A}}\right\} _{\mathcal{A}\in \mathrm{CA}},\;\;\psi _{%
\mathcal{A}}\in \frak{Z}^{2}\left[ \mathbf{A}_{1}\right] ,
\end{equation*}
such that $\omega _{\mathcal{A},\mathcal{B}}=\frak{j}\left( \psi _{\mathcal{B%
}}^{!},\psi _{\mathcal{A}},\mathbb{I}^{\otimes 2}\right) .\;\;\;\blacksquare
$
\end{definition}

Since Eq. $\left( \ref{m3}\right) $, such cochains $\omega _{\mathcal{A},%
\mathcal{B}}$ are in $\frak{Z}^{2}\left[ \mathbf{B}_{1}^{\ast }\otimes
\mathbf{A}_{1}\otimes \mathbf{B}_{1}\right] $. To give an example, in the
TTP case with $\widehat{\tau }_{\mathcal{A},\mathcal{B}}=id\otimes \sigma _{%
\mathcal{B}}^{!}\otimes \sigma _{\mathcal{A}}$, the cochain $\psi _{\mathcal{%
A}}$ would be given by the assignment
\begin{equation}
a_{k_{1}}...a_{k_{r}}\otimes a_{k_{r+1}}...a_{k_{r+s}}\mapsto
a_{k_{1}}...a_{k_{r}}\otimes \left( \sigma _{\mathcal{A}}^{-r}\right)
_{k_{r+1}}^{j_{1}}...\left( \sigma {}_{\mathcal{A}}^{-r}\right)
\,_{k_{r+s}}^{j_{s}}\,a_{j_{1}}...a_{j_{s}}.  \label{fia}
\end{equation}

From the injection
\begin{equation}
\frak{Z}^{2}\left[ \mathbf{B}_{1}\right] ^{!}\times \frak{Z}^{2}\left[
\mathbf{A}_{1}\right] \times \frak{Z}^{2}\left[ \mathbf{B}_{1}\right]
\hookrightarrow \frak{Z}^{2}\left[ \mathbf{B}_{1}^{\ast }\otimes \mathbf{A}%
_{1}\right] \times \frak{Z}^{2}\left[ \mathbf{B}_{1}\right] ,  \label{3inc}
\end{equation}
we shall also regard $\omega _{\mathcal{A},\mathcal{B}}$ as a cochain
belonging to the latter set, depending on our convenience.

\begin{theorem}
If a category $\Omega ^{\cdot }$ is associated to a factorizable collection,
then $\Omega ^{\cdot }$ is a sub-semigroupoid of $\frak{H}\mathrm{CA}^{\circ
}$.
\end{theorem}

\begin{proof}
Consider the quantum spaces $\mathcal{A}$, $\mathcal{B}$ and $\mathcal{C}$,
and diagrams $\left\langle \varphi ,\mathcal{H}\right\rangle \in \Omega ^{%
\mathcal{A},\mathcal{B}}$ and $\left\langle \psi ,\mathcal{G}\right\rangle
\in \Omega ^{\mathcal{B},\mathcal{C}}$, with associated linear spaces (\emph{%
via }the functor\emph{\ }$\frak{F}$)
\begin{equation*}
\mathbf{H}_{1}^{\varphi }=span\left[ h_{i}^{j}\right] _{i,j=1}^{\frak{n},%
\frak{m}}\;\;and\;\;\mathbf{G}_{1}^{\psi }=span\left[ g_{i}^{j}\right]
_{i,j=1}^{\frak{m},\frak{p}}.
\end{equation*}
We are denoting by $\dim \mathbf{A}_{1}=\frak{n}$, $\dim \mathbf{B}_{1}=%
\frak{m}$ and $\dim \mathbf{C}_{1}=\frak{p}$ the dimensions of the generator
spaces defining $\mathcal{A}$, $\mathcal{B}$ and $\mathcal{C}$,
respectively. We must show that $\left\langle \left( I_{H}\circ \chi \right)
\,\varphi ,\mathcal{H}\circ \mathcal{G}\right\rangle $ (see Eq. $\left( \ref
{ma}\right) $) is an object of $\Omega ^{\mathcal{A},\mathcal{C}}$, and that
the objects $\left\langle \ell _{\mathcal{A}},\mathcal{K}\right\rangle $ are
in $\Omega ^{\cdot }$. That means the quantum space $\frak{F}\left\langle
\left( I_{H}\circ \chi \right) \,\varphi ,\mathcal{H}\circ \mathcal{G}%
\right\rangle $, generated by
\begin{equation*}
span\left[ \sum_{j=1...\frak{m}}h_{i}^{j}\otimes g_{j}^{k}\right] _{i,k=1}^{%
\frak{n},\frak{p}}\subset \mathbf{H}_{1}^{\varphi }\otimes \mathbf{G}%
_{1}^{\psi }\subset \mathbf{H}^{\varphi }\circ \mathbf{G}^{\psi }\mathbf{,}
\end{equation*}
is such that $\left( I_{H}\circ \chi \right) \,\varphi $ defines an arrow $%
\mathcal{A}\rightarrow \left( \frak{F}\left\langle \left( I_{H}\circ \chi
\right) \,\varphi ,\mathcal{H}\circ \mathcal{G}\right\rangle \circ \mathcal{C%
}\right) _{\omega }$ in $\mathrm{CA}$. To this end, let us introduce some
notation.

Denote by $\mathbf{h}$ and $\mathbf{g}$ the matrices with entries $%
h_{i}^{j}\in \mathbf{H}_{1}$ and $g_{i}^{j}\in \mathbf{G}_{1}$, and by $%
\mathbf{a}$, $\mathbf{b}$ and $\mathbf{c}$ the vectors whose components are $%
a_{i}\in \mathbf{A}_{1}$, $b_{i}\in \mathbf{B}_{1}$ and $c_{i}\in \mathbf{C}%
_{1}$. Since $\left\langle \varphi ,\mathcal{H}\right\rangle $ and $%
\left\langle \psi ,\mathcal{G}\right\rangle $ are elements of $\Omega
^{\cdot }$, $\omega _{\mathcal{A},\mathcal{B}}$ and $\omega _{\mathcal{B},%
\mathcal{C}}$ defines cochains in $\frak{Z}^{2}\left[ \mathbf{H}%
_{1}^{\varphi }\right] \times \frak{Z}^{2}\left[ \mathbf{B}_{1}\right] $ and
$\frak{Z}^{2}\left[ \mathbf{G}_{1}^{\psi }\right] \times \frak{Z}^{2}\left[
\mathbf{C}_{1}\right] $ (see Eq. $\left( \ref{3inc}\right) $), respectively.
The latter are given by
\begin{equation}
\omega _{\mathcal{A},\mathcal{B}}\left( \mathbf{h}_{r,s}\otimes \mathbf{b}%
_{r,s}\right) =\left[ \psi _{\mathcal{A}}\right] _{r,s}\cdot \mathbf{h}%
_{r,s}\cdot \left[ \psi _{\mathcal{B}}^{-1}\right] _{r,s}\otimes \mathbf{b}%
_{r,s}  \label{ab}
\end{equation}
and
\begin{equation}
\omega _{\mathcal{B},\mathcal{C}}\left( \mathbf{g}_{r,s}\otimes \mathbf{c}%
_{r,s}\right) =\left[ \psi _{\mathcal{B}}\right] _{r,s}\cdot \mathbf{g}%
_{r,s}\cdot \left[ \psi _{\mathcal{C}}^{-1}\right] _{r,s}\otimes \mathbf{c}%
_{r,s}  \label{bc}
\end{equation}
where the symbols $\mathbf{h}_{r,s}=\mathbf{h}_{r}\otimes \mathbf{h}_{s}$
and $\left[ \psi _{\mathcal{A}}\right] _{r,s}\cdot \mathbf{h}_{r,s}\cdot %
\left[ \psi _{\mathcal{B}}^{-1}\right] _{r,s}$ denote elements of the form
\begin{equation*}
h_{i_{1}}^{j_{1}}...h_{i_{r}}^{j_{r}}\otimes
h_{k_{1}}^{l_{1}}...h_{k_{s}}^{l_{s}}\in \left( \mathbf{H}_{1}^{\varphi
}\right) ^{\otimes r}\otimes \left( \mathbf{H}_{1}^{\varphi }\right)
^{\otimes s}
\end{equation*}
and
\begin{equation*}
\left( \psi _{\mathcal{A}}\right)
_{n_{1}...n_{r},m_{1}...m_{s}}^{i_{1}...i_{r},k_{1}...k_{s}}%
\;h_{i_{1}}^{j_{1}}...h_{i_{r}}^{j_{r}}\otimes
h_{k_{1}}^{l_{1}}...h_{k_{s}}^{l_{s}}\;\left( \psi _{\mathcal{B}%
}^{-1}\right) _{j_{1}...j_{r},l_{1}...l_{s}}^{p_{1}...p_{r},q_{1}...q_{s}},
\end{equation*}
respectively. Now, $\left\langle \left( I_{H}\circ \chi \right) \,\varphi ,%
\mathcal{H}\circ \mathcal{G}\right\rangle \in \Omega ^{\mathcal{A},\mathcal{C%
}}$ if and only if $\left( I_{H}\circ \chi \right) \,\varphi $ defines the
mentioned arrow in $\mathrm{CA}$, with $\omega $ given by
\begin{equation*}
\omega _{\mathcal{A},\mathcal{C}}\left( \left( \mathbf{h}\overset{.}{\otimes
}\mathbf{g}\right) _{r,s}\otimes \mathbf{c}_{r,s}\right) =\left[ \psi _{%
\mathcal{A}}\right] _{r,s}\cdot \left( \mathbf{h}\overset{.}{\otimes }%
\mathbf{g}\right) _{r,s}\cdot \left[ \psi _{\mathcal{C}}^{-1}\right]
_{r,s}\otimes \mathbf{c}_{r,s}.
\end{equation*}
Here $\overset{.}{\otimes }$ is denoting matrix contraction between $\mathbf{%
h}$ and $\mathbf{g}$. It follows from straightforward calculations that, if $%
\omega _{\mathcal{A},\mathcal{C}}$ is well defined on $\frak{F}\left\langle
\left( I_{H}\circ \chi \right) \,\varphi ,\mathcal{H}\circ \mathcal{G}%
\right\rangle \circ \mathcal{C}$, then $\left( I_{H}\circ \chi \right)
\,\varphi $ is a homogeneous linear map $\mathbf{A}\rightarrow \mathbf{H}%
\otimes \mathbf{G}\otimes \mathbf{C}$ defining the wanted morphism. So, let
us first show that. To this end, extend $\omega _{\mathcal{A},\mathcal{C}}$
to $\mathbf{H}^{\varphi }\circ \mathbf{G}^{\psi }\circ \mathbf{C}$ by
putting $\omega _{\mathcal{A},\mathcal{C}}\left( \mathbf{h}_{r,s}\otimes
\mathbf{g}_{r,s}\otimes \mathbf{c}_{r,s}\right) $ equal to
\begin{equation*}
\begin{array}{l}
\left[ \psi _{\mathcal{A}}\right] _{r,s}\cdot \mathbf{h}_{r,s}\otimes
\mathbf{g}_{r,s}\cdot \left[ \psi _{\mathcal{C}}^{-1}\right] _{r,s}\otimes
\mathbf{c}_{r,s}= \\
\\
=\left[ \psi _{\mathcal{A}}\right] _{r,s}\cdot \mathbf{h}_{r,s}\cdot \left[
\psi _{\mathcal{B}}^{-1}\right] _{r,s}\overset{.}{\otimes }\left[ \psi _{%
\mathcal{B}}\right] _{r,s}\cdot \mathbf{g}_{r,s}\cdot \left[ \psi _{\mathcal{%
C}}^{-1}\right] _{r,s}\otimes \mathbf{c}_{r,s}.
\end{array}
\end{equation*}
From the last expression, and recalling that $\omega _{\mathcal{A},\mathcal{B%
}}$ and $\omega _{\mathcal{B},\mathcal{C}}$ (given in $\left( \ref{ab}%
\right) $ and $\left( \ref{bc}\right) $) are admissible, it follows that $%
\omega _{\mathcal{A},\mathcal{C}}$ is admissible for $\mathcal{H}^{\varphi
}\circ \mathcal{G}^{\psi }\circ \mathcal{C}$, and also for the subspace $%
\frak{F}\left\langle \left( I_{H}\circ \chi \right) \,\varphi ,\mathcal{H}%
\circ \mathcal{G}\right\rangle \circ \mathcal{C}$. Then, $\left\langle
\left( I_{H}\circ \chi \right) \,\varphi ,\mathcal{H}\circ \mathcal{G}%
\right\rangle \in \Omega ^{\mathcal{A},\mathcal{C}}$.

Finally, we have to show the units $\left\langle \ell _{\mathcal{A}},%
\mathcal{K}\right\rangle $\textsf{\ }are objects of the corresponding
categories $\Omega ^{\mathcal{A}}$. We know that $\ell _{\mathcal{A}}$ is a
homogeneous linear map such that $\ell _{\mathcal{A}}\left( a\right)
=e^{n}\otimes a$, if $a\in \mathbf{A}_{n}$. In particular, we can write $%
\ell _{\mathcal{A}}\left( a_{i}\right) =e\delta _{i}^{j}\otimes a_{j}$.
Then, the cochain $\omega _{\mathcal{A}}$ for $\left\langle \ell _{\mathcal{A%
}},\mathcal{K}\right\rangle $ is given by
\begin{equation*}
\omega _{\mathcal{A}}\left( \mathbf{k}_{r,s}\otimes \mathbf{a}_{r,s}\right) =%
\left[ \psi _{\mathcal{A}}\right] _{r,s}\cdot \mathbf{k}_{r,s}\cdot \left[
\psi _{\mathcal{A}}^{-1}\right] _{r,s}\otimes \mathbf{a}_{r,s}
\end{equation*}
with
\begin{equation*}
\mathbf{k}_{r,s}=e^{r}\,\delta _{i_{1}}^{j_{1}}...\delta
_{i_{r}}^{j_{r}}\otimes e^{s}\,\delta _{k_{1}}^{l_{1}}...\delta
_{k_{s}}^{l_{s}}.
\end{equation*}
Accordingly $\omega _{\mathcal{A}}$ for $\left\langle \ell _{\mathcal{A}},%
\mathcal{K}\right\rangle $ is the identity map, and $\left( \mathcal{K}\circ
\mathcal{A}\right) _{\omega }=\mathcal{K}\circ \mathcal{A}$.
\end{proof}

\bigskip

The following result is immediate.

\begin{proposition}
Let us call $\circ _{\Omega }$ the partial product associated to above
mentioned semigroupoid structure of $\Omega ^{\cdot }$. The embedding $\frak{%
P}^{\Omega }:\Omega ^{\cdot }\hookrightarrow \mathrm{CA}:\left\langle
\varphi ,\mathcal{H}\right\rangle \mapsto \mathcal{H}$ preserves the
respective \emph{(}partial\emph{)} products and units, i.e.
\begin{equation*}
\frak{P}^{\Omega }\,\circ _{\Omega }=\circ \,\left( \frak{P}^{\Omega }\times
\frak{P}^{\Omega }\right) \;\;\;and\;\;\;\frak{P}^{\Omega }\left\langle \ell
_{\mathcal{A}},\mathcal{K}\right\rangle =\mathcal{K}.\;\;\;\blacksquare
\end{equation*}
\end{proposition}

\section{The generalized coHom objects}

For $\Omega ^{\mathcal{A},\mathcal{B}}$ to have initial objects we need an
additional condition on $\omega _{\mathcal{A},\mathcal{B}}$.

\begin{theorem}
If $\Omega ^{\cdot }$ is associated to a factorizable collection given by $%
\left\{ \psi _{\mathcal{A}}\right\} _{\mathcal{A}\in \mathrm{CA}}$ such that
each $\frak{i}\psi _{\mathcal{A}}=\frak{i}\psi $ is \textbf{2nd} $\mathcal{A}
$-\textbf{admissible}, then each $\Omega ^{\mathcal{A},\mathcal{B}}$ have
initial object
\begin{equation*}
\underline{hom}^{\Omega }\left[ \mathcal{B},\mathcal{A}\right] =\mathcal{B}_{%
\frak{i}\psi }\triangleright \mathcal{A}_{\frak{i}\psi }=\left( \mathcal{B}%
\triangleright \mathcal{A}\right) _{\frak{j}\left( \frak{i}\psi ^{!},\frak{i}%
\psi \right) }.
\end{equation*}
In particular, $\underline{hom}^{\Omega }\left[ \mathcal{K},\mathcal{A}%
\right] =\mathcal{A}_{\frak{i}\psi }$ and $\underline{hom}^{\Omega }\left[
\mathcal{K},\mathcal{K}\right] =\mathcal{K}$; thus,
\begin{equation*}
\underline{hom}^{\Omega }\left[ \mathcal{B},\mathcal{A}\right] =\underline{%
hom}^{\Omega }\left[ \mathcal{K},\mathcal{B}\right] \triangleright
\underline{hom}^{\Omega }\left[ \mathcal{K},\mathcal{A}\right]
.\;\;\;\blacksquare
\end{equation*}
\end{theorem}

Before going to the proof, let us make some remarks. Since $\psi _{\mathcal{A%
}}\in \frak{Z}^{2}\left[ \mathbf{A}_{1}\right] $, there exists a primitive $%
\theta _{\mathcal{A}}\in \frak{P}^{1}\left[ \mathbf{A}_{1}\right] $ such
that $\psi _{\mathcal{A}}=\partial \theta _{\mathcal{A}}$, $\frak{i}\psi _{%
\mathcal{A}}=\partial \theta _{\mathcal{A}}^{-1}$ and $\psi _{\mathcal{A}%
}^{!}=\partial ^{!}\theta _{\mathcal{A}}^{!}$ (see \S \textbf{3.2.1},\textbf{%
\ }\S \textbf{3.3.1 }and \S \textbf{4.2.1}, respectively). In addition, if $%
\mathbf{I}=\bigoplus_{n\geq 2}\mathbf{I}_{n}$ is the graded ideal related to
$\mathcal{A}$, we have from \S \textbf{3.2.2} that (provided $\frak{i}\psi $
is admissible)
\begin{equation}
\mathbf{I}_{\frak{i}\psi ,n}=\theta \left( \mathbf{I}_{n}\right) ,\;\;%
\mathbf{I}_{\frak{i}\psi ,n}^{\perp }=\theta \left( \mathbf{I}_{n}\right)
^{\perp }=\theta ^{!}\left( \mathbf{I}_{n}^{\perp }\right) .  \label{m}
\end{equation}

\begin{proof}
\textbf{(of theorem)} We shall show $\mathcal{B}_{\frak{i}\psi
}\triangleright \mathcal{A}_{\frak{i}\psi }$ defines an object of $\Omega ^{%
\mathcal{A},\mathcal{B}}$ and then that is initial.

Let us note that, given $\psi ,\varphi \in \frak{Z}^{2}$, with $\psi
=\partial \theta $ and $\varphi =\partial \chi $, it follows that
\begin{eqnarray*}
\frak{j}\left( \frak{i}\psi ,\frak{i}\varphi \right)  &=&\frak{j}\left(
\partial \left( \theta ^{-1}\right) ,\partial \left( \chi ^{-1}\right)
\right) =\frak{j}\left( \partial \left( \theta ^{-1},\chi ^{-1}\right)
\right) = \\
&=&\partial \left( \frak{j}\left( \theta ,\chi \right) ^{-1}\right) =\frak{ij%
}\left( \psi ,\varphi \right) .
\end{eqnarray*}
Also recall, if $\frak{i}\psi $ is (2nd) $\mathcal{A}$-admissible, then $%
\psi $ is (2nd) $\mathcal{A}_{\frak{i}\psi }$-admissible (see \textbf{Prop.
14} of \S \textbf{3.3.1}).

By \textbf{Theor. 16} of \S \textbf{4.2.2}, using the 2nd $\mathcal{A}$%
-admissibility of $\frak{i}\psi _{\mathcal{A}}=\frak{i}\psi $,
\begin{equation*}
\begin{array}{l}
\left( \mathcal{B}_{\frak{i}\psi }\triangleright \mathcal{A}_{\frak{i}\psi
}\right) \circ \mathcal{B}=\left( \mathcal{B}\triangleright \mathcal{A}%
\right) _{\frak{j}\left( \frak{i}\psi ^{!},\frak{i}\psi \right) }\circ
\mathcal{B}=\left( \mathcal{B}\triangleright \mathcal{A}\right) _{\frak{ij}%
\left( \psi ^{!},\psi \right) }\circ \mathcal{B} \\
\\
=\left( \left( \mathcal{B}\triangleright \mathcal{A}\right) \circ \mathcal{B}%
\right) _{\frak{ij}\left( \psi ^{!},\psi ,\mathbb{I}^{\otimes 2}\right)
}=\left( \left( \mathcal{B}\triangleright \mathcal{A}\right) \circ \mathcal{B%
}\right) _{\frak{i}\omega },
\end{array}
\end{equation*}
because $\omega =\omega _{\mathcal{A},\mathcal{B}}=\frak{j}\left( \psi _{%
\mathcal{B}}^{!},\psi _{\mathcal{A}},\mathbb{I}^{\otimes 2}\right) $. Then,
\begin{equation*}
\left( \left( \mathcal{B}_{\frak{i}\psi }\triangleright \mathcal{A}_{\frak{i}%
\psi }\right) \circ \mathcal{B}\right) _{\omega }=\left( \left( \left(
\mathcal{B}\triangleright \mathcal{A}\right) \circ \mathcal{B}\right) _{%
\frak{i}\omega }\right) _{\omega }=\left( \mathcal{B}\triangleright \mathcal{%
A}\right) \circ \mathcal{B},
\end{equation*}
and consequently the map $\delta :\mathcal{A}\rightarrow \left( \mathcal{B}%
\triangleright \mathcal{A}\right) \circ \mathcal{B}:a_{i}\mapsto
z_{i}^{j}\otimes b_{j}$ (with $z_{i}^{j}=b^{j}\otimes a_{i}$), defining the
coevaluation of the proper coHom object $\underline{hom}\left[ \mathcal{B},%
\mathcal{A}\right] =\mathcal{B}\triangleright \mathcal{A}$, gives also a
morphism $\delta :\mathcal{A}\rightarrow \left( \left( \mathcal{B}_{\frak{i}%
\psi }\triangleright \mathcal{A}_{\frak{i}\psi }\right) \circ \mathcal{B}%
\right) _{\omega }$. Moreover, since the equality $\frak{F}\left\langle
\delta ,\mathcal{B}_{\frak{i}\psi }\triangleright \mathcal{A}_{\frak{i}\psi
}\right\rangle =\mathcal{B}_{\frak{i}\psi }\triangleright \mathcal{A}_{\frak{%
i}\psi }$, the pair $\left\langle \delta ,\mathcal{B}_{\frak{i}\psi
}\triangleright \mathcal{A}_{\frak{i}\psi }\right\rangle $ is in $\Omega ^{%
\mathcal{A},\mathcal{B}}$. Let us show such a pair is an initial object.

Suppose $\mathcal{A}$ and $\mathcal{B}$ have related dimensions $\dim
\mathbf{A}_{1}=\frak{n}$ and $\dim \mathbf{B}_{1}=\frak{m}$, and ideals $%
\mathbf{I}$ and $\mathbf{J}$, respectively. Let us consider the vector space
\begin{equation*}
\mathbf{D}_{1}\doteq span\left[ z_{i}^{j}\right] _{i,j=1}^{\frak{n},\frak{m}}
\end{equation*}
and the linear map $\delta _{1}:a_{i}\mapsto z_{i}^{j}\otimes b_{j}$. Under
the identification $b^{j}\otimes a_{i}=z_{i}^{j}$, the cochain $\omega _{%
\mathcal{A},\mathcal{B}}=\omega $ defines a counital 2-cocycle in $\frak{C}%
^{2}\left[ \mathbf{D}_{1}\otimes \mathbf{B}_{1}\right] $, and $\delta _{1}$
can be extended to an algebra homomorphism
\begin{equation*}
\delta _{1}^{\otimes }:\mathbf{A}_{1}^{\otimes }\rightarrow \left( \mathbf{D}%
_{1}^{\otimes }\circ \mathbf{B}_{1}^{\otimes }\right) _{\omega }.
\end{equation*}
Analogous calculations that enable us to arrive at Eq. $\left( 3.4\right) $
of \cite{g}, show that
\begin{equation*}
\delta _{1}^{\otimes }\left( \mathbf{a}_{r}\right) =\left[ \theta _{\mathcal{%
A}}\right] _{r}\cdot \mathbf{z}_{r}\cdot \left[ \theta _{\mathcal{B}}^{-1}%
\right] _{r}\overset{.}{\otimes }\mathbf{b}_{r},
\end{equation*}
(where we are using notation of previous theorem), or in coordinates,
\begin{equation*}
\delta _{1}^{\otimes }\left( a_{i_{1}}...a_{i_{r}}\right) =\left( \theta _{%
\mathcal{A}}\right)
_{i_{1}...i_{r}}^{j_{1}...j_{r}}\;z_{j_{1}}^{l_{1}}...z_{j_{r}}^{l_{r}}\;%
\left( \theta _{\mathcal{B}}^{-1}\right)
_{l_{1}...l_{r}}^{k_{1}...k_{r}}\otimes b_{k_{1}}...b_{k_{s}}.
\end{equation*}
In fact, identifying $\mathbf{D}_{1}^{\otimes }$ with $\left[ \mathbf{B}%
_{1}^{\ast }\otimes \mathbf{A}_{1}\right] ^{\widehat{\otimes }}$, we have
\begin{equation*}
\left[ \theta _{\mathcal{A}}\right] _{r}\cdot \mathbf{z}_{r}\cdot \left[
\theta _{\mathcal{B}}^{-1}\right] _{r}=\theta _{\mathcal{B}}^{!}\left(
\mathbf{b}_{r}^{\ast }\right) \otimes \theta _{\mathcal{A}}\left( \mathbf{a}%
_{r}\right) ;\;\;\mathbf{b}_{r}^{\ast }=b^{k_{1}}...b^{k_{r}}.
\end{equation*}
In this notation, the generators of $\mathbf{I}$, $\mathbf{J}$ and $\mathbf{J%
}^{\perp }$ can be written
\begin{equation*}
\left[ R_{\lambda _{r}}\right] _{r}\cdot \mathbf{a}_{r},\;\;\;\left[ S_{\mu
_{r}}\right] _{r}\cdot \mathbf{b}_{r}\;\;\;and\;\;\;\mathbf{b}_{r}^{\ast
}\cdot \left[ S^{\perp \omega _{r}}\right] _{r},
\end{equation*}
respectively. To have a well defined map $\delta $ from $\delta
_{1}^{\otimes }$ when $\mathbf{A}_{1}^{\otimes }$ and $\mathbf{B}%
_{1}^{\otimes }$ are quotient by their corresponding ideals, $\mathbf{D}%
_{1}^{\otimes }$ must be quotient by the elements
\begin{equation}
\left[ R_{\lambda _{r}}\right] _{r}\cdot \left[ \theta _{\mathcal{A}}\right]
_{r}\cdot \mathbf{z}_{r}\cdot \left[ \theta _{\mathcal{B}}^{-1}\right]
_{r}\cdot \left[ S^{\perp \omega _{r}}\right] _{r}.  \label{idz}
\end{equation}
Now, it is clear that a pair $\left\langle \varphi ,\mathcal{H}\right\rangle
\in \frak{H}\mathrm{CA}^{\circ }$ is in $\Omega ^{\mathcal{A},\mathcal{B}}$
if and only if there exist elements $h_{i}^{j}\in \mathbf{H}_{1}\subset
\frak{F}\left\langle \varphi ,\mathcal{H}\right\rangle $ satisfying
relations $\left( \ref{idz}\right) $ (replacing $\mathbf{z}_{r}$ by $\mathbf{%
h}_{r}$). But the elements of $\left( \ref{idz}\right) $ span precisely the
space (from Eq. $\left( \ref{m}\right) $)
\begin{equation*}
\theta _{\mathcal{B}}^{!}\left( \mathbf{J}_{r}^{\perp }\right) \otimes
\theta _{\mathcal{A}}\left( \mathbf{I}_{r}\right) =\mathbf{J}_{\frak{i}\psi
,r}^{\perp }\otimes \mathbf{I}_{\frak{i}\psi ,r},
\end{equation*}
which generates algebraically the ideal related to $\mathcal{B}_{\frak{i}%
\psi }\triangleright \mathcal{A}_{\frak{i}\psi }$. Then, the function $%
z_{i}^{j}\mapsto h_{i}^{j}$ can be extended to an arrow $\mathcal{B}_{\frak{i%
}\psi }\triangleright \mathcal{A}_{\frak{i}\psi }\rightarrow \frak{F}%
\left\langle \varphi ,\mathcal{H}\right\rangle $. But this is the unique
arrow in $\Omega ^{\mathcal{A},\mathcal{B}}$ that can be defined between
these objects, that is to say, $\left\langle \delta ,\mathcal{B}_{\frak{i}%
\psi }\triangleright \mathcal{A}_{\frak{i}\psi }\right\rangle $ is an
initial object of $\Omega ^{\mathcal{A},\mathcal{B}}$.

Finally, since the cochains of $\frak{C}^{\bullet }\left[ \Bbbk \right] $
are always (2nd) $\mathcal{K}$-admissible, in particular the primitive ones $%
\frak{P}^{1}\left[ \Bbbk \right] $, it follows that any twisting of $%
\mathcal{K}$ is isomorphic to $\mathcal{K}$. Then, the last claim of the
theorem follows immediately from the first one.
\end{proof}

\bigskip

Note that the 2nd admissibility condition for each $\frak{i}\psi _{\mathcal{A%
}}$ replaces the automorphism property of the related $\sigma _{\mathcal{A}}$
of the TTP case. Such a condition is immediate in the TTP case, because the
properties defining a twisting map imply the related 2-cocycles are
anti-bicharacters (c.f. \textbf{\S 2.2.2} and \textbf{\S 4.1.2}).

Now, a couple of immediate corollaries.

\begin{corollary}
The gauge equivalence \emph{(}see \textbf{\S 3.3}\emph{)} $\underline{hom}%
^{\Omega }\left[ \mathcal{B},\mathcal{A}\right] \backsim \underline{hom}%
\left[ \mathcal{B},\mathcal{A}\right] $ is valid for all conic quantum
spaces $\mathcal{B},\mathcal{A}$. \ \ \ $\blacksquare $
\end{corollary}

\begin{corollary}
In the category of quadratic quantum spaces $\mathrm{QA}$ \emph{(}and any $%
\mathrm{CA}^{m}$\emph{)}, the initial objects of $\Omega ^{\mathcal{A},%
\mathcal{B}}$ are isomorphic to
\begin{equation*}
\underline{hom}^{\Omega }\left[ \mathcal{B},\mathcal{A}\right] =\left(
\mathcal{B}_{\frak{i}\psi }\right) ^{!}\bullet \mathcal{A}_{\frak{i}\psi
}=\left( \mathcal{B}^{!}\right) _{\frak{i}\psi ^{!}}\bullet \mathcal{A}_{%
\frak{i}\psi }=\left( \mathcal{B}^{!}\bullet \mathcal{A}\right) _{\frak{j}%
\left( \frak{i}\psi ^{!},\frak{i}\psi \right) }.\;\;\;\blacksquare
\end{equation*}
\end{corollary}

Under the hypothesis mentioned in previous theorem, and from the
semigroupoid structure of $\Omega ^{\cdot }$ compatible with the monoid in $%
\mathrm{CA}$, we have the annunciated result:

\begin{theorem}
The assignment $\left( \mathcal{B},\mathcal{A}\right) \mapsto \underline{hom}%
^{\Omega }\left[ \mathcal{B},\mathcal{A}\right] $ define an $\mathrm{CA}$%
-cobased category with arrows
\begin{equation*}
\underline{hom}^{\Omega }\left[ \mathcal{C},\mathcal{A}\right] \rightarrow
\underline{hom}^{\Omega }\left[ \mathcal{B},\mathcal{A}\right] \circ
\underline{hom}^{\Omega }\left[ \mathcal{C},\mathcal{B}\right] ,
\end{equation*}
the cocomposition, and for $\underline{end}^{\Omega }\left[ \mathcal{A}%
\right] $ the counit epimorphism
\begin{equation*}
\underline{end}^{\Omega }\left[ \mathcal{A}\right] \twoheadrightarrow
\mathcal{K}\;\;\;/\;\;\;z_{i}^{j}\mapsto \delta _{i}^{j}\;e,
\end{equation*}
and the monomorphic comultiplication
\begin{equation*}
\underline{end}^{\Omega }\left[ \mathcal{A}\right] \hookrightarrow
\underline{end}^{\Omega }\left[ \mathcal{A}\right] \circ \underline{end}%
^{\Omega }\left[ \mathcal{A}\right] \;\;\;/\;\;\;z_{i}^{j}\mapsto
z_{i}^{k}\otimes z_{k}^{j}.\;\;\;\blacksquare
\end{equation*}
\end{theorem}

\end{document}